\documentclass[12pt,twoside]{article}
\usepackage{amsmath}
\usepackage{amsfonts}
\usepackage{amsthm}
\usepackage{amssymb}
\usepackage{graphicx}
\usepackage{multicol}

\begin{document}
\title{{\normalsize{\bf Free ribbon lemma for surface-link}}}
\author{{\footnotesize Akio Kawauchi}\\
{\footnotesize{\it Osaka Central Advanced Mathematical Institute, Osaka Metropolitan University}}\\
{\footnotesize{\it Sugimoto, Sumiyoshi-ku, Osaka 558-8585, Japan}}\\
{\footnotesize{\it kawauchi@omu.ac.jp}}}
\date\,
\maketitle
\vspace{0.25in}
\baselineskip=10pt
\newtheorem{Theorem}{Theorem}[section]
\newtheorem{Conjecture}[Theorem]{Conjecture}
\newtheorem{Lemma}[Theorem]{Lemma}
\newtheorem{Sublemma}[Theorem]{Sublemma}
\newtheorem{Proposition}[Theorem]{Proposition}
\newtheorem{Corollary}[Theorem]{Corollary}
\newtheorem{Claim}[Theorem]{Claim}
\newtheorem{Definition}[Theorem]{Definition}
\newtheorem{Example}[Theorem]{Example}

\begin{abstract} A free surface-link is a surface-link whose fundamental group is a free group 
not necessarily meridian-based. 
Free ribbon lemma says that every free sphere-link in the 4-sphere is a ribbon sphere-link. Four different proofs of Free ribbon lemma are explained. 
The first proof is done in an earlier paper. The second proof is done by showing that 
there is an O2-handle basis of a ribbon surface-link.
The third proof is done by removing the commutator relations from a Wirtinger presentation of a free 
group, which a paper on another proof of Free ribbon lemma complements. 
The fourth proof is given by the special case of the proof of the result that every free surface-link is a ribbon surface-link which is a stabilization of a free ribbon sphere-link. 
As a consequence, it is shown that a surface-link is a sublink of a free surface-link if and only if 
it is a stabilization of a ribbon sphere-link. 

\phantom{x}

\noindent{\it Keywords}: Free ribbon lemma,\, Free surface-link,\, Ribbon sphere-link,\, Stabilization.

\noindent{\it Mathematics Subject Classification 2010}: Primary 57N13; Secondary 57Q45

\end{abstract}

\baselineskip=15pt

\bigskip

\noindent{\bf 1. Introduction}

A {\it surface-link} is a closed, possibly disconnected, oriented surface $F$ smoothly embedded in the 4-sphere $S^4$, and it is called a {\it surface-knot} if $F$ is connected.
If $F$ consists of 2-spheres $F_i\,(i=1,2,\dots,r)$, then $F$ is called 
a {\it sphere-link} (or an $S^2$-{\it link}) of $r$ components. 
It is shown that a surface-link $F$ is a trivial surface-link (i.e., bounds disjoint handlebodies in $S^4$) if the fundamental group $\pi_1(S^4\setminus F,x_0)$ is a meridian-based free group, \cite{K21}, \cite{K23-1}, \cite{K23-2}. 
A surface-link $F$ is {\it ribbon} if $F$ is obtained from a trivial $S^2$-link $O$ in $S^4$ 
by surgery along a smoothly embedded disjoint 1-handle system $h^O$ on $O$, 
\cite{KSS}, \cite{Yajima62}, \cite{Yajima64}, \cite{Yana}. 
A surface-link $F$ in the 4-sphere $S^4$ is {\it free} if the fundamental group 
$\pi_1(S^4\setminus F,x_0)$ is a (not necessarily meridian-based) free group. 
In this paper, four different proofs of the following {\it Free ribbon lemma} and its 
generalization to a general free surface-link are explained.

\phantom{x}

\noindent{\bf Free ribbon lemma.} Every free $S^2$-link in $S^4$ is a ribbon $S^2$-link.

\phantom{x}

Free ribbon lemma leads to the following conjectures: 
Poincar{\'e} conjecture, \cite{K24-2}, \cite{Per}, \cite{P1}. \cite{P2}. 
J. H. C. Whitehead asphericity conjecture for aspherical 2-complex, \cite{Hw}, \cite{K24-1}, 
\cite{K24-4}, \cite{Wh}. 
Kervaire conjecture on group weight, \cite{GR}, \cite{K24-3}, \cite{Ker}, \cite{MKS}, \cite{Kly}. 
The first proof is given \cite{K24-1}. For convenience, an outline of the first proof is explained here.

\phantom{x}

\noindent{\bf First proof of Free ribbon lemma.} 
Let $L_i\,(i=1,2,\dots,r)$ be the components of a free $S^2$-link $L$ in $S^4$. 
By a base change of the free fundamental group $\pi_1(S^4\setminus L, x_0)$, 
take a basis $x_i\,(i=1,2,\dots,r)$ of $\pi_1(S^4\setminus L, x_0)$ inducing 
a meridian basis of $L$ in $H_1(S^4\setminus L;Z)$, \cite{MKS}. 
Let $Y$ be the 4-manifold obtained from $S^4$ by surgery along $L$, which is diffeomorphic to 
the connected sum of $r$ copies $S^1\times S^3_i\, (i=1,2,\dots,r)$ of $S^1\times S^3$, 
\cite{K23-3}, \cite{K24-1}. 
Under a canonical isomorphism $\pi_1(S^4\setminus L,x_0)\to \pi_1(Y,x_0)$, 
the factors $S^1\times p_i\, (i=1,2,\dots,r)$ of $S^1\times S^3_i\, (i=1,2,\dots,r)$ 
with suitable paths to the base point $x_0$ represent the basis $x_i\, (i=1,2,\dots,r)$.
Let $k_i\,(i=1,2,\dots,r)$ be the loop system in $Y$ produced from the components 
$L_i\,(i=1,2,\dots,r)$ by the surgery. By using the fact that any homotopy deformations of 
$k_i\,(i=1,2,\dots,r)$ in $Y$ do not change the link type of the surface-link $L$ in $S^4$, 
the loop system $k_i\,(i=1,2,\dots,r)$ is homotopically deformed in $Y$ so that 
the surface-link $L$ in $S^4$ obtained from the deformed loop system $k_i\,(i=1,2,\dots,r)$ 
by back surgery is a ribbon surface-link in $S^4$, completing the proof of Free ribbon lemma.

\phantom{x}

To explain the second and third proofs of Free ribbon lemma, the notion of an O2-handle basis of a surface-link 
is needed, \cite{K21}, \cite{K24-6}. An {\it O2-handle pair}  on a surface-link $F$ in $S^4$ is a pair 
$(D\times I, D' \times I)$ of 2-handles 
$D\times I$, $D' \times I$ on $F$ in $S^4$ which intersect orthogonally only with the attaching parts 
$(\partial D)\times I$, $(\partial D')\times I$ to $F$, so that 
the intersection $Q=(\partial D)\times I\cap (\partial D')\times I$ is a square. 
Let $(D\times I, D'\times I)$ be an O2-handle pair on a surface-link $F$. 
Let $F(D\times I)$ and $F(D'\times I)$ be the surface-links obtained from $F$ by the surgeries along $D\times I$ and $D'\times I$, respectively. 
Let $F(D\times I, D'\times I)$ be the surface-link which is 
the union $\delta\cup F^c_{\delta}$ of the plumbed disk 
\[\delta=\delta_{D\times I,D'\times I}
=D\times \partial I\cup Q\cup D'\times \partial I\]
and the surface 
$F^c_{\delta}= \mbox{cl}(F\setminus (\partial D\times I\cup \partial D'\times I))$. 
The surface-links $F(D\times I), F(D'\times I)$ and $F(D\times I, D'\times I)$ are equivalent 
surface-links, \cite{K21}. 
An {\it O2-handle basis} of a surface-link $F$ is a disjoint system of O2-handle pairs 
$(D_i\times I, D'_i\times I)\, (i=1,2,\dots,r)$ on $F$ in $S^4$ such that the boundary loop 
pair system $(\partial D_i, \partial D'_i)\, (i=1,2,\dots,r)$ of the core disk system $(D_i, D'_i)\, (i=1,2,\dots,r)$ of $(D_i\times I, D'_i\times I)\, (i=1,2,\dots,r)$ is a spin loop basis for $F$ in $S^4$, which is a system of a spin loop basis of every component $F_i$ of $F$. Note that 
there is a spin loop basis for every surface-knot in $F$, \cite{HiK}.
In this paper, for simplicity, an O2-handle basis 
$(D_i\times I, D'_i\times I)\, (i=1,2,\dots,r)$ for $F$ is denoted by $(D\times I, D'\times I)$.  
The surgery surface-link of $F$ by $(D_i\times I, D'_i\times I)\, (i=1,2,\dots,r)$ is denoted by 
$F(D\times I, D'\times I)$. 
The following theorem is shown for the second and third proofs of Free ribbon lemma.

\phantom{x}

\noindent{\bf Theorem~1.1.} For every free ribbon surface-link $F$ in $S^4$, there is an O2-handle basis 
$(D\times I,D'\times I)$ on $F$ in $S^4$ such that $D\times I$ belongs to the 1-handle 
system of the ribbon surface-link $F$. 

\phantom{x}

The second proof of Free ribbon lemma is explained as follows.

\phantom{x}

\noindent{\bf Second proof of Free ribbon lemma.} 
Let $L$ be a free $S^2$-link. Then there is a ribbon surface-link $F$ such that 
the fundamental group $\pi_1(S^4\setminus F,x_0)$ is isomorphic to the free fundamental group 
$\pi_1(S^4\setminus L,x_0)$ by a meridian-preserving isomorphism, \cite{K24-5}. By Theorem~1.1, the surgery surface-link $L'=F(D\times I, D'\times I)$ is a ribbon $S^2$-link, \cite{K21}, \cite{K24-6}. Then there is a meridian-preserving isomorphism 
$\pi_1(S^4\setminus L',x_0)\to \pi_1(S^4\setminus L,x_0)$ on free groups, which implies that 
$L'$ is equivalent to $L$, \cite{K24-1}, \cite{K24-5}. 
Thus, $L$ is a ribbon $S^2$-link, completing the proof of Free ribbon lemma.

\phantom{x}

The third proof of Free ribbon lemma is related to a Wirtinger presentation of a free group.
A finite group presentation 
$(x_1, x_2, \dots, x_n |\, R_1, R_2,\dots, R_m)$ 
is a Wirtinger presentation if $R_j=W_j x_{s_j} {W_j}^{-1} x_{t_j}^{-1}$ forf some word $W_j$ in the free group 
$(x_1, x_2, \dots, x_n)$ and  some indexes $s_j, t_j$ in 
$\{1,2,\dots,n\}$ for every $j\,(j=1,2,\dots,m)$. The relator $R_j$ is a {\it commutator relation} if 
$x_{s_j}=x_{t_j}$. It is well-known that a Wirtinger presentation of a finitely presented group $G$ 
with $H_1(G;Z)\cong Z^r$ is always equivalent (without changing the generating set) to a Wirtinger presentation $P$ such that 
the Wirtinger presentation $P'$ obtained by removing all the commutator relations from $P$ 
has deficiency $r$. Such a Wirtinger presentation $P$ is called a {\it normal} Wirtinger presentation.
The following corollary is obtained from Theorem~1.1. 

\phantom{x}

\noindent{\bf Corollary~1.2.} If a free group $G$ of rank $r$ has a normal Wirtinger presentation 
$P$, then $G$ has the Wirtinger presentation $P'$ of deficiency $r$ obtained from $P$ 
by removing all the commutator relations. 

\phantom{x}

\noindent{\bf Proof of Corollary~1.2 assuming Theorem~1.1.}  For a free 
group $G$ of rank $r$, 
let $P=(x_1, x_2, \dots, x_n|\, R_1, R_2,\dots, R_m)$ be a normal Wirtinger presentation of  
$G$ such that the relators $R_j\, (n-r+1\leq j \leq m)$ are the commutator relations.
Let $O$ be a trivial $S^2$-link of $n$ components in $S^4$ such that the 
meridian basis of the free fundamental group $\pi_1(S^4\setminus O,x_0)$ are
identified with $x_i\, (i=1,2,\dots,n)$. Let $h_j\, (1\leq j \leq m)$ be the 1-handles on $O$ 
indicated by the relators $R_j\, (1\leq j \leq m)$. By the van Kampen theorem, 
the ribbon surface-link $F$ in $S^4$ obtained by 
surgery along $h_j\, (1\leq j \leq m)$ has the normal Wirtinger presentation $P$ 
of the fundamental group $\pi_1(S^4\setminus F,x_0)$ with the meridian generators 
set $\{x_1, x_2, \dots, x_n\}$, \cite{K07}, \cite{K08}. 
Let $L$ be the ribbon surface-link obtained from $O$ by surgery along 
the 1-handles $h_j\, (1\leq j \leq n-r)$, which is a ribbon $S^2$-link of $r$ components. 
The fundamental group $\pi_1(S^4\setminus L,x_0)$ has the Wirtinger presentation $P'$ of deficiency $r$ obtained from $P$ by removing all the commutator relations. 
By Theorem~1.1, the 1-handles $h_j\, (n-r+1\leq j \leq m)$ on $L$ are trivial 1-handles, so that 
$\pi_1(S^4\setminus L,x_0)$ is isomorphic to $\pi_1(S^4\setminus F,x_0)$ by a meridian-preserving isomorphism. This completes the proof of Corollary~1.2 assuming Theorem~1.1. 

\phantom{x}

The author has published a paper on another proof of Free ribbon lemma, which this paper complements, \cite{K24-5}. The third proof of Free ribbon lemma is nothing but the proof of the paper
except for adding to it the assertion of Corollary~1.5 which was missing from it. 
For convenience, an outline of the third proof is explained here.

\phantom{x}

\noindent{\bf Third proof of Free ribbon lemma.} 
Let $L$ be a free $S^2$-link of $r$ components. 
Since the fundamental group $G= \pi_1(S^4\setminus L, x_0)$
is a free group with $H_1(G;Z)=Z^r$ and $H_2(G;Z)=0$, there is a normal Wirtinger presentation 
$P$ of $G$ whose generator set comes from meridians of $L$ in $S^4$,  
\cite{K24-5}, \cite{Yajima70}. Note that there is also another method to find such a normal Wirtinger presentation $P$ using a normal form of $L$ in $S^4$, \cite{Kam01},  \cite{K07}, \cite{K08}, \cite{KSS}. 
Let $L'$ be a ribbon $S^2$-link given by the Wirtinger presentation $P'$ obtained from $P$ 
by removing all the commutators. By Corollary~1.2, there is a meridian-preserving isomorphism 
$\pi_1(S^4\setminus L',x_0)\to \pi_1(S^4\setminus L,x_0)$, so that $L'$ is equivalent to 
$L$. Thus, $L$ is a ribbon $S^2$-link, completing the proof of Free ribbon lemma.

\phantom{x}

The fourth proof of Free ribbon lemma is given by a direct proof of the following theorem.

\phantom{x}

\noindent{\bf Theorem~1.3.} Every free surface-link $F$ in $S^4$ is a ribbon surface-link in $S^4$. 

\phantom{x} 

\noindent{\bf Fourth proof of Free ribbon lemma.} It is obtained 
from Theorem~1.3  by restricting $F$ to every 
free $S^2$-link, completing the proof of Free ribbon lemma.

\phantom{x} 

Thus, after the proofs of Theorems~1.1 and 1.3, there are four different proofs of Free ribbon lemma. To generalize the free ribbon lemma to a free surface-link, the notion of a stabilization of 
a surface-link is needed, \cite{K21}, \cite{K24-6}. 
A {\it stabilization} of a surface-link $L$ is a connected sum 
$ F= L\# _{k=1}^s T_k$ of $L$ and a system of trivial torus-knots 
$T_k\, (k=1,2,\dots, s)$. 
By granting $s=0$, a surface-link $L$ itself is regarded as a stabilization 
of $L$. 
Free ribbon lemma is generalized to a general free surface-link as follows.

\phantom{x}

\noindent{\bf Corollary~1.4.} Every free surface-link $F$ in $S^4$ is a stabilization of a free ribbon $S^2$-link $L$ in $S^4$. 

\phantom{x}

\noindent{\bf Proof of Corollary~1.4 assuming Theorems~1.1 and 1.3.}
Theorems~1.1 and 1.3 imply that every free surface-link $F$ is a ribbon surface-link and a stabilization of a free $S^2$-link $L$, \cite{K21}. By Free ribbon lemma, the free $S^2$-link $L$ is a ribbon $S^2$-link. 
This result also follows directly, \cite{K24-6}. 
This completes the proof of Corollary~1.4 assuming Theorems~1.1 and 1.3.

\phantom{x}

It is shown that an $S^2$-link $L$ is a sublink of a free $S^2$-link if and 
only if 
$L$ is a ribbon $S^2$-link, \cite{K24-1}. The following corollary generalizes this property to a general surface-link. 

\phantom{x}

\noindent{\bf Corollary~1.5.} A surface-link $L$ in $S^4$ is a sublink of a free surface-link $F$ in $S^4$ if and only if $L$ is a stabilization of a ribbon $S^2$-link in $S^4$. 

\phantom{x}

\noindent{\bf Proof of Corollary~1.5 assuming Corollary~1.4.} If $L$ is a sublink of a free surface-link $F$, then $L$ is a stabilization of a ribbon $S^2$-link since every free surface-link 
is a stabilization of a free ribbon $S^2$-link by Corollary~1.4. Conversely, 
if $L$ is a stabilization of a ribbon $S^2$-link, then $L$ is a sublink of a stabilization of a free ribbon 
$S^2$-link which is a free surface-link $F$ since every ribbon $S^2$-link is a sublink of a free $S^2$-link. This completes the proof of Corollary~1.5 assuming Corollary~1.4.

\phantom{x}

\noindent{\bf 2. Proofs of Theorems~1.1 and 1.3.} 

Let $F$ be a {\it free} surface-link in $S^4$ with components 
$F_i\, (i=1,2,\dots, r)$. 
Let $N(F)=\cup_{i=1}^r N(F_i)$ be a tubular neighborhood of 
$F=\cup_{i=1}^r F_i$ in $S^4$ which is a trivial normal disk bundle $F\times D^2$ over $F$, where 
$D^2$ denotes the unit disk of complex numbers of norm $\leqq 1$. 
Let $E=E(F)=\mbox{cl}(S^4\setminus N(F))$ be the exterior of $F$ in $S^4$. 
The boundary 
$\partial E=\partial N(F)=\cup_{i=1}^r \partial N(F_i)$ 
of the exterior $E$ is a trivial normal circle bundle over $F=\cup_{i=1}^r F_i$. 
Identify $\partial N(F_i)=F_i\times S^1$ for $S^1=\partial D^2$ such that 
the composite inclusion 
\[F_i\times 1\to \partial N(F_i) \to \mbox{cl}(S^4\setminus N(F_i))\]
induces the zero-map in the integral first homology. 
The following lemma uses the assumption that the fundamental group $\pi_1(E,x_0)$ is a free group of rank $r$ and the fact that the first homology group $H_1(E;Z)$ is a free abelian group of rank $r$ with meridian basis.

\phantom{x}

\noindent{\bf Lemma~2.1.} 
The composite inclusion $F_i\times 1\to\partial N(F_i) \to E$ is null-homotopic for all $i$. 

\phantom{x}

\noindent{\bf Proof of Lemma~2.1.} Since $\partial N(F_i)=F_i\times S^1$, 
the fundamental group elements between the factors $F_i\times 1$ and $q_i\times S^1$ are commutive. Let $a_i\, (i=1,2,\dots,r)$ be embedded edges with common vertex $x_0$ in $E$ 
such that $a_i\setminus\{x_0\}\, (i=1,2,\dots,r)$ are mutually disjoint and 
$a_i\cap(\cup_{j=1}^r F_j\times 1)=p_i\times 1$ for a point $p_i$ of $F_i\times 1$. 
The surface $F_i\times 1$ in $\partial N(F_i)=F_i\times S^1$ is chosen so that 
the inclusion $F_i\times 1\to\mbox{cl}(S^4\setminus N(F_i))$ 
induces the zero-map in the integral first homology. 
Since $H_1(E;Z)$ is a free abelian group of rank $r$ with meridian basis and 
$\pi_1(E,x_0)$ is a free group of rank $r$, the image of the homomorphism 
$\pi_1(a_i\cup F_i\times S^1,x_0)\to \pi_1(E,x_0)$ is an infinite cyclic group generated by the 
homotopy class $[a_i\cup p_i \times S^1]$. 
This implies that the inclusion $F_i\times 1\to E$ is null-homotopic. 
This completes the proof of Lemma~2.1. 
$\square$ 

\phantom{x}

By using the free group $\pi_1(E,x_0)$ of rank $r$, let 
$\Gamma=\cup_{i=1}^r a_i \cup C_i$ 
be a connected graph in the interior $\mbox{Int}(E)$ of $E$ 
consisting of embedded edges $a_i\, (i=1,2,\dots,r)$ with the common base point $x_0$ and disjoint embedded circles 
$C_i\, (i=1,2,\dots,r)$ such that 

\medskip

\noindent{(1)} the half-open edges $a_i\setminus \{x_0\} \, (i=1,2,\dots,r)$ are mutually disjoint and 
$a_i \cap (\cup_{j=1}^r C_j)=v_i$, a point in $C_i$ for every $i$, 

\medskip

\noindent{(2)} the inclusion $i: (\Gamma,x_0)\to (E,x_0)$ induces an isomorphism 
$i_{\#}: \pi_1(\Gamma,x_0)\to \pi_1(E,x_0)$, and

\medskip

\noindent{(3)} the homology class $[p_i\times S^1]= [C_i]$ in $H_1(E;Z)$ for all $i$. 

\phantom{x}

In fact, by (2), the homotopy classes $[a_i\cup C_i]\, (i=1,2,\dots,r)$ form a basis of the free 
group $\pi_1(E,q_0)$. (3) is obtained by a base change of the free group $\pi_1(E,x_0)$, \cite{MKS}. 
Since $\Gamma$ is a $K(\pi,1)$-space, there is a piecewise-linear map 
$f:(E,x_0)\to (\Gamma,x_0)$ inducing the inverse isomorphism 
$f_{\#}=(i_{\#})^{-1}: \pi_1(E,q_0)\to \pi_1(\Gamma,q_0)$ of $i_{\#}$, and by the homotopy extension property, 
the restriction of $f$ to $\Gamma$ is the identity map, \cite{Spanier}. 
By Lemma~2.1, the restriction of $f$ to $\partial E$ is homotopic to the composite map 
\[g: \partial E=F\times S^1 \to \cup_{i=1}^r q_i\times S^1 \to \Gamma\]
such that the first map $F\times S^1 \to \cup_{i=1}^r q_i\times S^1$ 
is induced from the constant map $F_i \to \{q_i\}$ for all $i$ 
and the second map $\cup_{i=1}^r q_i\times S^1 \to \Gamma$ is defined by the map $f$. 
By using a boundary collar of $\partial E$ in $E$, assume that the piecewise-linear map 
$f:(E,x_0)\to (\Gamma,x_0)$ defines the map $g: \partial E \to \Gamma$. 
For a non-vertex point $p_i$ of $C_i$, the preimage $V_i=(f)^{-1}(p_i)$ is a 
bi-collard compact oriented proper piecewise-linear 3-manifold in $E$. 
Take the compact 4-manifold $E'$  obtained from $E$ by splitting along $V=\cup_{i=1}^r V_i$  to be connected. Then  join the components in  each 
$V_i$ with 1-handles in $E'$.  By these modifications, $V_i$ is assumed to be connected for all i within a homotopic deformation of  $f$. 
The boundary $\partial V_i$ is the disjoint union 
$P_i(F)$ of $m_{ij}$ parallel copies $m_{ij} F_j$
of $F_j\times 1$ for all $j\,(j=1,2,\dots,r)$ in $S^4$, where $m_{ii}$ is an odd integer and $m_{ij}$ for distinct $i,j$ is an even integer. 
An {\it anti-parallel surface-link} of $F_j$ in $S^4$ 
is the boundary  
$\partial(F_j\times J)$ of a normal line bundle $F_j\times J$ of 
$F_j$ in  $S^4$ such that the natural homology homomorphism 
$H_1(F\times 1;Z)\to H_1(S^4\setminus F\times 0;Z)$ is the zero-map, where $J=[0,1]$. 
Note that the boundary  $\partial (F^{(0)}_j\times J)$ for a 
compact once-punctured surface $F^{(0)}_j$ of $F_j$  is a trivial surface-knot 
in $S^4$ since $F^{(0)}_j\times J$ is a handlebody. 
Let $n_{ii}=(m_{ii}-1)/2$ and $n_{ij}=m_{ij}/2$ for distinct $i,j$.    
The surface-link $P_i(F)=\partial V_i$ consists of a component identified with $F_i$ and some anti-parallel 
surface-links  of $F_j$ in $P_i(F)$, denoted by 
$P_{ijk}(F)\, (j=1,2,\dots, r; k=1,2,\dots, n_{ij})$.
Let $P(F)=\cup_{i=1}^r P_i(F)$ be the surface-link in $S^4$.

The proof of Theorem~1.1 is done as follows.

\phantom{x}

\noindent{\bf Proof of Theorem~1.1.} 
A ribbon surface-link $F$ is obtained from a trivial $S^2$-link $O$ in 
$S^4$ by surgery along a disjoint 1-handle system $h^O$ on $O$. Then,
the surface-link $P(F)$ of a free ribbon surface-link $F$ is a ribbon 
surface-link obtained from a trivial $S^2$-link $P(O)$, a parallel link of $O$,  in 
$S^4$ by surgery along a disjoint 1-handle system $P(h^O)$, a parallel 1-handle system of $h^O$,   on $P(O)$. 
Let $(\ell,\ell')$ be a spin simple loop basis  for $P(F)$ consisting a spin simple loop basis of every component of $P(F)$, where the spin loop system $\ell$ is taken to be the boundary of a  disk system 
$D$ in  a transverse disk system of the 1-handle system $P(h^O)$ and the spin loop system $\ell'$ 
is  the union of an arc system parallel to the core arc system of $P(h^O)$ and an arc system 
in $P(F)$. Because of the isomorphism $f_{\#}$,  the compact 4-manifold $E'$ 
is simply connected, so that  every loop of $\ell'$ is null-homotopic in $E'$ and hence 
bounds an immersed disk system in $E'$. 
Let $(\ell_1,\ell'_1)$ be a loop pair in $(\ell,\ell')$ meeting at a boundary point of a 1-handle $h^O_1$  in $P(h^O)$, and $\delta_1'$ an immersed disk in $E'$ with $\partial \delta_1'=\ell'_1$. 
By thickening the arc $\ell'_1\cap P(h^O)$ to the 1-handle $h^O_1$, 
the  loop pair $(\ell_1,\ell'_1)$ bounds an immersed O2-handle pair $(D_1\times I, \delta_1'\times I)$ on  $P(F)$ and hence there is  an O2-handle pair $(D_1\times I, D_1'\times I)$ on $P(F)$ 
such that the boundary loop pair $(\partial D_1, \partial D'_1)$ of the core disk pair $(D_1, D'_1)$ is  
the  loop pair $(\ell_1,\ell'_1)$  by Recovery Lemma, \cite{K23-1}.
By using the 3-ball $D_1\times I\cup D_1'\times I$, every loop of $\ell'\setminus {\ell'_1}$ bounds 
an immersed disk system in $E'$ not meeting $D_1'\times I$. Thus, by continuing this procedure, 
the spin simple loop basis $(\ell,\ell')$ bounds an O2-handle basis $(D\times I, D'\times I)$ on $P(F)$  
such that the 2-handle system $D\times $I belongs to $P(h^O)$, whose subsystem on $F$ gives a desired
 O2-handle basis on $F$ in $S^4$. This completes the proof of Theorem~1.1.

\phantom{x}

Let $h_i$ be a disjoint 1-handle system on $P_i(F)$ embedded in $V_i$ 
such that the surface $P_i(F; h_i)$ obtained from $P_i(F)$ by surgery along $h_i$ is connected and the genus of $P_i(F; h_i)$ is equal to the total genus of $P_i(F)$. 
Assume that a just one 1-handle $h^F_i$ of $h_i$ attaches to $F_i$. 
Let $h=\cup_{i=1}^r h_i$ be the 1-handle system  on $P(F)$, and 
 $P(F; h)=\cup_{i=1}^r P_i(F; h_i)$ a surface-link in $S^4$. 
By further taking a disjoint 1-handle system 
$h'_i$ on $P_i(F; h_i)$ embedded in $V_i$, the closed surface $P_i(F; h_i, h'_i)$ obtained from 
$P_i(F; h_i)$ by surgery along $h'_i$ bounds a handlebody $V'_i$ in $V_i$, so that the surface-link 
$P(F; h, h')=\cup_{i=1}^r P_i(F; h_i, h'_i)$ is a trivial surface-link in $S^4$. 
Since the compact 4-manifold $E'$ 
is simply connected,  the 1-handle system 
$h'=\cup _{i=1}^r h'_i$ is a trivial 1-handle system on the surface-link 
$P(F; h)$ in $S^4$, \cite{HoK}, \cite{K23-2}. 
Thus, the surface-link $P(F; h)$ is a trivial surface-link in $S^4$, \cite{K21}, \cite{K23-1}. 
The proof of Theorem~1.3 is done as follows.

\phantom{x}

\noindent{\bf Proof of Theorem~1.3.} 
Since $E'$ is simply connected, the 1-handle system $h_i$ 
can be  chosen so that  a trivial surface-knot  $F_{ijk}$ is obtained from the  anti-parallel surface-link $P_{ijk}(F)$ in  $P_i(F)$  by surgery along 
a unique 1-handle $h_{ijk}$ in $h_i$.  
Let $h(0)$ be the system of the 1-handles  
$h_{ijk}\, (i, j=1,2,\dots, r;  k=1,2,\dots, n_{ij})$ in $h$,  
and $h(1)$  the complementary system of $h(0)$ in $h$. 
Let $P(F;h(0))$ be the surface-link obtained from $P(F)$ by surgery 
along $h(0)$, which consists  of the surface-link $F$ and  the trivial 
surface-knots  $F_{ijk}\, (k=1,2,\dots, n_{ij}; i, j=1,2,\dots, r)$.   
The trivial surface-link $P(F;h)$ is obtained from 
the surface-link $P(F;h(0))$ by surgery along  the 1-handle system 
$h(1)$. As stated just before the proof of Theorem~1.3, 
the surface-link $P(F; h, h')$ bounds a  handlebody system 
$V'$ in $V$. 
Note that  $h'$ is a trivial 1-handle system on $P(F; h)$. Then there is a disjoint handlebody system $U$ in $S^4$ with $\partial U=P(F; h)$ extending the handlebody system $V'$ by adding  a 2-handle system 
$e\times I$ which makes an O2-handle system together with a thickened transverse disk system $m(h')$ of $h'$, \cite{K21}, \cite{K23-1}.  
Let $d$ be a  transverse  disk system of the 1-handle system $h(1)$. 
In general, the disk system $d$ meets the core disk 
system $e$ of the 2-handle system $e\times I$ transversely in finite points in $S^4$, but 
the handlebody system $U$ is isotopically deformed so to have 
$d\cap U =\partial d$ 
by an isotopic deformation of $U$ in a neighborhood of the 3-disk system 
$e\times I\cup m(h')\times I$ in $S^4$, \cite{K23-1}. 
For $h(1)=d\times I$, the union $U\cup h(1)$ is a compact oriented 
3-manifold with boundary $P(F;h(0))$ obtained from $U$ by  adding the 
2-handle system $h(1)$ since 
$h(1)\cap U=h(1)\cap \partial U = (\partial d)\times I$.  
Let $(\ell,\ell')$ be a spin loop basis of $P(F;h)$ given for $P(F;h(0))$ 
such that when restricted to every trivial surface-knot $F_{ijk}$ 
it becomes a standard  spin loop basis.  
Since $P(F;h)$ is a trivial surface-link, there is an orientation-preserving diffeomorphism $w$ of $S^4$ sending $P(F;h)$ to the boundary $\partial W$  
of a standard handlebody system $W$ in $S^4$ such that 
the spin loop basis $(w(\ell),w(\ell'))$ of $\partial W$ 
is  a meridian-longitude pair system of $W$, \cite{Hiro}, \cite{K21}. 
Let $U(W)=w^{-1}(W)$ be the handlebody system in $S^4$ with 
$\partial U(W)=P(F;h)$.  
The spin loop basis $(w(\ell),w(\ell'))$ of $\partial W$ 
bounds a core disk-pair system $(\delta,\delta')$ of an O2-handle basis 
$(\delta\times I, \delta'\times I)$ of the trivial surface-link $\partial W$ in $S^4$, where $\delta$ denotes a meridian disk system of $W$. 
Thus, the spin loop basis $(\ell,\ell')$ of $P(F; h)$ bounds the core disk pair system $(D,D')$ of the O2-handle basis 
$(D\times I,D'\times I)=(w^{-1}(\delta)\times I, w^{-1}(\delta')\times I)$ 
on $P(F; h)$ in $S^4$ with $D$ in $U(W)$, so that 
$d\cap D=h(1)\cap D=\emptyset$. 
The handlebody system $U(W)$ is isotopically 
deformed to be $U(W)=U$ in $S^4$.
Consider  the 3-manifold $U'\cup h(1)$ obtained from 
the 3-manifold $U\cup h(1)$ by splitting along the disk system $D$, which 
is a multi-punctured 3-sphere system. 
Since $\partial(U\cup h(1))=P(F;h(0))$, the boundary $\partial (U'\cup h(1))$ is obtained 
from  $P(F;h(0))$ by surgery along the 
2-handle system $D\times I$ and consists of the $S^2$-link $S$ 
(obtained from $F$) and the  $S^2$-knots  
$S_{ijk}$ (obtained from $F_{ijk}$)  for all $i, j, k$. 
The  $S^2$-knots $S_{ijk}$ for all $i, j, k$ are shown to form  a trivial 
$S^2$-link in $S^4$. 
To see this,  consider a collar $V\times J$ of $V$  in $S^4$ to 
move (via a boundary collar of $V$) the anti-parallel surface 
$P_{ijk}(F)$  and the 1-handle $h_{ijk}$ into a  disjoint subcollar 
$V\times J_{ijk}$ of  $V\times J$  for each $i, j, k$.  
Then by forgetting the 1-handle system $h(1)$, the trivial surface-knots $F_{ijk}$ for all $i, j, k$ (not containing $F$) in the boundary of $U\cup h(1)$ are isotopic  to the boundaries of disjoint handlebodies $U_{ijk}$ for all $i, j, k$ in $S^4$.   
For the transverse disk system $d$ of  
the 1-handle system $h(1)$ on  $P(F;h(0))$, 
there is  a proper disk system  $d_U$ in $U$ with $\partial d_U=\partial d$ 
which splits $U$ into disjoint handlebodies $U^0$ and 
$U^0_{ijk}$ for all $i, j, k$  such that  the boundaries $\partial U^0$ and 
$\partial U^0_{ijk}$ contain the multi-punctured surfaces of 
$F$ and $F_{ijk}$, respectively  and $U^0_{ijk}$ contains the subsystem 
$D_{ijk}$ for  $F_{ijk}$ of the 2-handle core disk system $D$ on $P(F;h(0))$.  
Deform the handlebody $U^0_{ijk}$ into  the handlebody $U_{ijk}$ for all $i, j, k$, 
so that  the $S^2$-knots $S_{ijk}$ for all $i, j, k$ are isotopic to 
the boundaries  of  the disjoint 3-balls 
$U'_{ijk}$  for all $i, j, k$ obtained from the disjoint handlebodies $U_{ijk}$  
for all $i, j, k$ by splitting along  the disk systems $D_{ijk}$ for all $i, j, k$. 
Thus, the system of the $S^2$-knots $S_{ijk}$ for all $i, j, k$ forms  a trivial $S^2$-link in $S^4$, as desired.  
The multi-punctured 3-sphere system $U'\cup h(1)$ then means that  the $S^2$-link  $S$ is a ribbon $S^2$-link in $S^4$, \cite{K24-1}.  
Because $F$ is obtained from the ribbon $S^2$-link $S$ by surgery along a 1-handle system, 
the surface-link $F$ is a ribbon surface-link in $S^4$. This completes the proof of Theorem~1.3. 

\phantom{x}

\noindent{\bf Acknowledgements.} 
The author has tried to find a simple algebraic method to prove Corollary~1.2  but failed. He is grateful to Igor Mineyev for asking questions to this algebraic method. This work was partly supported by JSPS KAKENHI Grant Number JP21H00978 and  MEXT Promotion of Distinctive Joint Research Center Program JPMXP0723833165 and Osaka Metropolitan University Strategic Research Promotion Project (Development of International Research Hubs). 

\phantom{x}


\begin{thebibliography}{99}

\bibitem{GR} F. Gonz{\'a}lez Acu{\~n}a and A. Ram{\'i}rez, A knot-theoretic equivalent of the Kervaire conjecture, J. Knot Theory Ramifications 15 (2006), 471-478. 

\bibitem{Hiro} S. Hirose, On diffeomorphisms over surfaces trivially 
embedded in the 4-sphere, Algebraic and Geometric Topology 2 (2002), 791-824. 

\bibitem{HiK} J. A. Hillman and A. Kawauchi, Unknotting orientable surfaces in the 4-sphere, J. Knot Theory Ramifications 4 (1995), 213-224.

\bibitem{HoK} F. Hosokawa and A. Kawauchi, Proposals for unknotted 
surfaces in four-space, Osaka J. Math. 16 (1979), 233-248. 

\bibitem{Hw} J. Howie, Some remarks on a problem of J. H. C. Whitehead, Topology 22 (1983), 475-485.

\bibitem{Kam01} S. Kamada, Wirtinger presentations for higher dimensional manifold knots obtained from diagrams, Fund. Math. 168(2001), 105-112. 


\bibitem{K07} A. Kawauchi, On the surface-link groups, Intelligence of low dimensional topology 2006, Series on knots and everything 40 (2007), 157-164, World Sci. publ.

\bibitem{K08} A. Kawauchi, The first Alexander Z[Z]-modules of 
surface-links and of virtual links, Geometry \& Topology 
Monographs 14 (2008), 353-371. 

\bibitem{K18} A. Kawauchi, Faithful equivalence of equivalent ribbon 
surface-links, Journal of Knot Theory and Its Ramifications  27 (2018), 1843003 (23 pages). 

\bibitem{K21} A. Kawauchi, Ribbonness of a stable-ribbon surface-link, I. A stably trivial surface-link, Topology and its Applications 301 (2021), 107522 (16pages).

\bibitem{K23-1} A. Kawauchi, Uniqueness of an orthogonal 2-handle pair on a surface-link, Contemporary Mathematics (UWP) 4 (2023), 182-188. 

\bibitem{K23-2} A. Kawauchi, Triviality of a surface-link with meridian-based free fundamental group, Transnational Journal of Mathematical Analysis and Applications 11 (2023), 19-27. 

\bibitem{K23-3} A. Kawauchi, Smooth homotopy 4-sphere, WSEAS Transactions on Mathematics 22 (2023), 690-701.

\bibitem{K24-1} A. Kawauchi, Ribbonness of Kervaire's sphere-link in homotopy 4-sphere and its consequences to 2-complexes, J Math Techniques Comput Math 3(4) (2024), 01-08 (online). 

\bibitem{K24-2} A. Kawauchi, Classical Poincar{\'e} conjecture via 4D topology, J Math Techniques Comput Math  3(4) (2024), 1-7 (online). 

\bibitem{K24-3} A. Kawauchi, Kervaire conjecture on weight of group via fundamental group of ribbon sphere-link, J Math Techniques Comput Math 3(4) (2024), 1-3 (online).

\bibitem{K24-4} A. Kawauchi, Whitehead aspherical conjecture via ribbon sphere-link, J Math Techniques Comput Math 3(5) (2024), 01-10 (online).

\bibitem{K24-5} A. Kawauchi, Another proof of free ribbon lemma, J Math Techniques Comput Math 3(9) (2024), 01-03 (online).

\bibitem{K24-6} A. Kawauchi, Ribbonness of a stable-ribbon surface-link, II. General case  (MDPI) Mathematics 13 (3), 402 (2025), 1-11.  
https://doi.org/10.3390/math13030402. 

\bibitem{KSS} A. Kawauchi, T. Shibuya, S. Suzuki, Descriptions on surfaces 
in four-space I : Normal forms, Mathematics Seminar Notes, Kobe University 10 (1982), 75-125; 
II: Singularities and cross-sectional links, 
Mathematics Seminar Notes, Kobe University 11 (1983), 31-69.
https://sites.google.com/view/kawauchiwriting. 

\bibitem{Ker} M. A. Kervaire, On higher dimensional knots, Differential and combinatorial topology 27 (1965), 105-119, Princeton Univ. Press. 

\bibitem{Kly} Ant. A. Klyachko, A funny property of a sphere and equations over groups, Comm. Algebra 21 (1993), 2555-2575.

\bibitem{MKS} W. Magnus, A. Karrass and D. Solitar, Combinatorial group theory: Presentations of groups in terms of generators and relations, Interscience Publishers (1966). 

\bibitem{Per} G. Perelman, Ricci flow with surgery on three-manifolds. arXiv: math.
DG/0303109 v1, 10 Mar 2003.

\bibitem{P1} H. Poincar{\'e}, Second compl{\'e}ment {\`a} l'Analysis Sitis, 
Proc. London Math. Soc. 32 (1900), 277-308. 

\bibitem{P2} H. Poincar{\'e}, Cinqui{\`e}me compl{\'e}ment {\`a} l'Analysis Sitis, Rend. Circ. Mat. Palermo 18 (1904), 45-110. 

\bibitem{Spanier} E. H. Spanier, Algebraic topology, MacGraw Hill (1966).

\bibitem{Wh} J. H. C. Whitehead, On adding relations to homotopy groups, Ann. Math. 42 (1941), 409-428.

\bibitem{Yajima62} T. Yajima, On the fundamental groups of knotted 2-manifolds in the 4-space, 
J. Math. Osaka City Univ.  13 (1962), 63-71. 

\bibitem{Yajima64} T. Yajima, On simply knotted spheres in $R^4$, Osaka J. Math.  1 (1964), 133-152.


\bibitem{Yajima70} T. Yajima, Wirtinger presentations of knot groups, Proc. Japan Acad. 46 (1970), 997-1000.

\bibitem{Yana} T. Yanagawa, On ribbon 2-knots; the 3-manifold bounded by 
the 2-knot, Osaka J. Math. 6 (1969), 447-164.

\end{thebibliography}
\end{document}